# Measurement of returns to scale with weight restrictions: How to deal with the occurrence of multiple supporting hyperplanes?

**Mahmood Mehdiloozad**

*Department of Mathematics, College of Sciences, Shiraz University, Shiraz 71454, Iran*
(E-mail: m.mehdiloozad@gmail.com)

**Kaoru Tone**

*National Graduate Institute for Policy Studies, Tokyo, Japan*
(E-mail: tone@grips.ac.jp)

**Biresh K. Sahoo**

*Xavier Institute of Management, Xavier University, Bhubaneswar 751013, India*
(E-mail: biresh@ximb.ac.in)

(✉) **Corresponding author**: M. Mehdiloozad

Ph.D. Candidate
Department of Mathematics
College of Sciences
Shiraz University
Golestan Street | Adabiat Crossroad | Shiraz 71454 | Iran

## Abstract

While measuring returns to scale in data envelopment analysis (DEA), the occurrence of multiple supporting hyperplanes has been perceived as a crucial issue. To deal effectively with this in weigh restrictions (WR) framework, we first precisely identify the two potential sources of its origin in the non-radial DEA setting. If the firm under evaluation *P* is WR-efficient, the non-full-dimensionality of its corresponding *P-face*—a face of minimum dimension that contains *P*—is the unique source of origin (problem Type I). Otherwise, the occurrence of multiple WR-projections or, correspondingly, multiple *P*-faces becomes the other additional source of origin (problem Type II). To the best of our knowledge, while problem Type I has been correctly addressed in the literature, the simultaneous occurrences of problems Types I and II have not effectively been coped with. Motivated by this, we first show that problem Type II can be circumvented by using a *P*-face containing all the *P*-faces. Based on this finding, we then devise a two-stage linear programming based procedure by extending a recently developed methodology by [Mehdiloozad, M., Mirdehghan, S. M., Sahoo, B. K., & Roshdi, I. (2015). On the identification of the global reference set in data envelopment analysis. *European Journal of Operational Research*, *245*, 779–788]. Our proposed method inherits all the advantages of the recently developed method and is computationally efficient. The practical applicability of our proposed method is demonstrated through a real-world data set of 80 Iranian secondary schools.

## 1 Introduction

One of the most important aspects in applied production analysis of firms (decision making units or DMUs) is the measurement of *returns to scale* (RTS) since its informational contents can provide important insights to firm managers making operational decisions in strengthening their competitive position (Tone & Sahoo, 2003; Sahoo & Tone, 2003, 2013, 2015). The economic concept of RTS was firstly introduced by Banker (1984) and Banker, Charnes, and Cooper (1984) into the nonparametric framework of data envelopment analysis (DEA)—a performance measurement method pioneered by Charnes, Cooper, and Rhodes (1978, 1979). Since then, it was extensively explored in the literature from both theoretical and practical aspects. See, e.g., Atici and Podinovski (2012), Banker and Thrall (1992), Banker, Bardhan, and Cooper (1996), Banker, Chang, and Cooper (1996), Banker, Cooper, Seiford, Thrall, and Zhu (2004), Førsund (1996), Fukuyama (2000), Golany and Yu (1997), Krivonozhko, Førsund, and Lychev (2012, 2014), Mehdiloozad et al. (2015), Podinovski and Førsund (2010), Podinovski, Førsund, and Krivonozhko (2009), Sahoo and Tone (2015), Sahoo, Kerstens, and Tone (2012), Sahoo, Zhu, and Tone (2014), Sahoo, Zhu, Tone, and Klemen (2014), Sueyoshi and Sekitani (2007a, 2007b), Tone (1996, 2005), Tone and Sahoo (2004, 2005, 2006), Zarepisheh and Soleimani-damaneh (2008), and Zarepisheh, Soleimani-damaneh, and Pourkarimi (2006), among others.

As is well-known, the BCC (Banker et al., 1984) model imposes no restriction on the weights attached to the inputs/outputs of firms except the non-negativity conditions. Because of this total weight flexibility, two serious drawbacks are associated with the use of the BCC model to actual situations: (1) the resulting optimal weights often take on zero values for the inputs and outputs of firms lying on the horizontal and vertical facets of input and output isoquants and (2) in other cases big differences in magnitudes of weights may exist. For more details, readers may refer to Ali, Cook, and Seiford (1991), Pedraja-Chaparro, Salinas-Jiminez, & Smith (1997), Roll and Golany (1993), Roll, Cook, & Golany (1991), and Thompson, Langemeier, Lee, C. T., Lee, E., and Thrall (1990). To address the above-mentioned drawbacks, the *weight restriction* (WR) research area has been introduced in the DEA literature to explore the incorporation of additional constraints on the magnitudes of weights. See, e.g., Dyson and Thanassoulis (1988), Estellita Lins, Moreira da Silva, and Lovell (2007), Podinovski (2001, 2004a, 2004b, 2007, 2013), Sarrico and Dyson (2004), Tone (2001), Tracy and Chen (2005) and Wong and Beasley (1990), among others.

While the measurement of RTS has been extensively discussed in the standard DEA environment, there exist a few research works studying this topic in the WR framework (see, e.g., Hosseinzadeh Lotfi, Jahanshahloo, and Esmaeili (2007), Korhonen, Soleimani-damaneh, and Wallenius (2011, 2013), Soleimani-damaneh, Jahanshahloo, Mehrabian, and Hasannasab (2010) and Tone (2001)). Tone (2001) attempted to measure RTS in the WR environment (hereafter referred to as WR-RTS). He introduced the notion of WR-RTS for a WR-efficient DMU with an unambiguous meaning and defined the WR-RTS for a WR-inefficient DMU as that of its WR-efficient WR-projection. Then, he extended the innovative method of Banker and Thrall (1992) for the measurement of the right- and left-hand WR-RTSs.

As in the conventional framework, the type and magnitude of the WR-RTS of a WR-efficient activity $P$ can be determined through the position(s) of the supporting hyperplane(s) binding at $P$. Obviously, there would be no problem if the hyperplane supporting at $P$ is unique. The difficulty, however, arises in the presence of multiplicity. To deal effectively with the multiplicity issue, the first step is to accurately identify its potential sources of origin. In this regard, we introduce the concept of *P-face* and define it as the face of minimum dimension that contains a WR-efficient activity $P$. Then, the intersection of technology set with all the supporting hyperplanes at $P$ is equal to the *P*-face, and the WR-RTS for $P$ mainly depends on the dimension of its associated *P*-face. More precisely, if the *P*-face is a '*Full Dimensional Efficient Facet*' (Olesen & Petersen, 1996, 2003), then a unique hyperplane is supporting at $P$ and the WR-RTS can be measured with no problem. Otherwise, the non-full-dimensionality of the *P*-face (hereafter referred to as *problem Type I*) raises the multiplicity issue as the first source.

It is thus clear that problem Type I becomes the unique source of the origin of the multiplicity issue when the DMU under evaluation is WR-efficient. However, while estimating the WR-RTS of a WR-inefficient DMU based on its WR-projection, the non-radial nature of the second phase of the WR model of Tone (2001) may become the other source. Precisely, in this case, multiple WR-projections and, correspondingly, multiple *P*-faces may be produced. Since each *P*-face is associated with the set of supporting hyperplanes at *P*, multiple sets of supporting hyperplanes may exist. Hereafter, this type of multiplicity is referred to as *problem Type II*.

To the best of our knowledge, while Tone's (2001) method correctly addresses problem Type I, there is no approach that effectively copes with the simultaneous occurrences of problems Types I and II. In the standard DEA environment, some studies (Krivonozhko et al., 2014; Mehdiloozad et al., 2015; Sueyoshi & Sekitani, 2007b) have recently succeeded to overcome the simultaneous occurrences of problems Types I and II. While the methods proposed in these studies are all interesting and theoretically correct, the method by Mehdiloozad et al. (2015) is computationally more efficient than the two others. Their proposed approach is a linear programming (LP) based method and contains two stages. In the first stage, the *global reference set* of the DMU under evaluation is identified through an LP model. A strict convex combination of all the identified reference units is then used as a projection for the measurement of RTS. The second stage involves the use of Banker et al.'s (2004) method for estimating the right- and left-hand RTSs of the obtained projection.

In this contribution, we extend Mehdiloozad et al.'s (2015) method to the WR framework. In this regard, we first introduce three types of WR reference set: *WR-unary reference set* (WR-URS), *WR-maximal reference set* (WR-MRS) and *WR-global reference set* (WR-GRS). For a given WR-projection *P*, the WR-URS is defined as the set of WR-efficient DMUs that are active in a specific convex combination producing *P*. The WR-MRS is the union of all the WR-URSs associated with *P*, and the WR-GRS is the union of WR-MRSs associated with all the WR-projections.

We show that a sufficient condition for circumventing problem Type II is the existence of a WR-projection $P^{max}$ such that $P^{max}$-face contains $P$-face for any WR-projection $P$. By a theorem, we also prove that $P^{max}$ always exists and, indeed, is a strict convex combination of the units in the WR-GRS. Then, we formulate an LP model for its identification. Finally, we devise a two-stage procedure to deal effectively with the simultaneous occurrences of problems Types I and II. In the first stage, we overcome problem Type II by finding $P^{max}$ via our proposed LP model. To resolve problem Type I, the WR-RTS of $P^{max}$ is then measured by using Tone's (2001) method at the second stage. We note that our proposed method inherits all the advantages of Mehdiloozad et al.'s (2015) and is computationally efficient.

The reminder of this paper is organized as follows. Section 2 deals with a brief review of the WR model and holds a discussion on measuring the WR-RTS. Section 3 first introduces the three concepts of the WR-URS, WR-MRS and WR-GRS; then, develops an LP problem for the identification of the WR-GRS; and finally, presents a two-stage procedure for characterizing the WR-RTS of inefficient DMUs. Section 4 applies our proposed method to an empirical data set. Section 5 concludes with some remarks.

## 2 Background

Throughout this paper, we deal with $n$ observed DMUs; each uses $m$ inputs to produce $s$ outputs. Let $J$ be the index set of the observed DMUs, i.e., $J = \{1,...,n\}$. For any $j \in J$, the input and output vectors of DMU$_j$ are denoted by $\mathbf{x}_j = (x_{1j},...,x_{mj})^T \in \mathbb{R}^m_{\geq 0}$ and $\mathbf{y}_j = (y_{1j},...,y_{sj})^T \in \mathbb{R}^s_{\geq 0}$, respectively. Further, $\mathbf{X} = [\mathbf{x}_1 \ ... \ \mathbf{x}_n]$ and $\mathbf{Y} = [\mathbf{y}_1 \ ... \ \mathbf{y}_n]$ are the input and output matrices, respectively. We denote vectors and matrices in bold letters, vectors in lower case and matrices in upper case. All vectors are column vectors. We denote by a subscript $T$ the transpose of vectors and matrices. We also use $\mathbf{0}_n$ and $\mathbf{1}_n$ to show $n$-dimensional vectors with the values of 0 and 1 in every entry, respectively.

## 2.1 The WR model

Assume that DMU$_o$ is assessed with respect to the technology set ($T$) defined as

$$T = \left\{ (\mathbf{x}, \mathbf{y}) \in \mathbb{R}_{\geq 0}^{m} \times \mathbb{R}_{\geq 0}^{s} \,\middle|\, \mathbf{x} \text{ can produce } \mathbf{y} \right\}. \tag{1}$$

Under variable returns to scale (VRS) and weight restriction assumptions, the non-parametric DEA-based representation of $T$ is set up as follows (Tone, 2001):

$$T_V^{WR} = \left\{ (\mathbf{x}, \mathbf{y}) \,\middle|\, \mathbf{X}\boldsymbol{\lambda} - \mathbf{P}\boldsymbol{\pi} \leq \mathbf{x},\ \mathbf{Y}\boldsymbol{\lambda} + \mathbf{Q}\boldsymbol{\tau} \geq \mathbf{y},\ \mathbf{1}_n^T \boldsymbol{\lambda} = 1,\ \boldsymbol{\lambda} \geq \mathbf{0}_n,\ \boldsymbol{\pi} \geq \mathbf{0}_k,\ \boldsymbol{\tau} \geq \mathbf{0}_l \right\}. \tag{2}$$

Here, $\mathbf{P}_{m \times k}$ and $\mathbf{Q}_{s \times l}$ are, respectively, associated with weight restrictions on inputs and outputs, and $\boldsymbol{\lambda}$ is called the *intensity vector*.

In reference to $T_V^{WR}$, the envelopment form of the *weight restriction (WR)* model (Tone, 2001) is represented as

$$\begin{aligned}
\min\quad & \theta - \varepsilon \left( \mathbf{1}_m^T \mathbf{s}^- + \mathbf{1}_s^T \mathbf{s}^+ \right) \\
\text{subject to}\quad & \\
& \mathbf{X}\boldsymbol{\lambda} - \mathbf{P}\boldsymbol{\pi} + \mathbf{s}^- = \theta \mathbf{x}_o, \\
& \mathbf{Y}\boldsymbol{\lambda} + \mathbf{Q}\boldsymbol{\tau} - \mathbf{s}^+ = \mathbf{y}_o, \\
& \mathbf{1}_n^T \boldsymbol{\lambda} = 1, \\
& \boldsymbol{\lambda} \geq \mathbf{0}_n,\ \boldsymbol{\pi} \geq \mathbf{0}_k,\ \boldsymbol{\tau} \geq \mathbf{0}_l, \\
& \mathbf{s}^- \geq \mathbf{0}_m,\ \mathbf{s}^+ \geq \mathbf{0}_s,
\end{aligned} \tag{3}$$

where $\varepsilon$ is a non-Archimedean infinitesimal. Let $\left( \theta^*, \boldsymbol{\lambda}^*, \boldsymbol{\pi}^*, \boldsymbol{\tau}^*, \mathbf{s}^{-*}, \mathbf{s}^{+*} \right)$ be an optimal solution to model (3). Then, *WR-efficiency* and *WR-improvement* are defined as follows.

**Definition 2.1** (*WR-efficiency*) DMU$_o$ is said to be WR-efficient if and only if $\theta^* = 1$, $\mathbf{s}^{-*} = \mathbf{0}_m$ and $\mathbf{s}^{+*} = \mathbf{0}_s$. Otherwise, it is called WR-inefficient.

**Definition 2.2** (*WR-improvement*) For a WR-inefficient DMU$_o$, a WR-projection is given by

$$P_o = (\hat{\mathbf{x}}_o, \hat{\mathbf{y}}_o) := (\mathbf{X}\boldsymbol{\lambda}^*, \mathbf{Y}\boldsymbol{\lambda}^*) = (\theta^* \mathbf{x}_o + \mathbf{P}\boldsymbol{\pi}^* - \mathbf{s}^{-*}, \mathbf{y}_o - \mathbf{Q}\boldsymbol{\tau}^* + \mathbf{s}^{+*}). \qquad (4)$$

**Theorem 2.1** (Tone, 2001) The activity $(\hat{\mathbf{x}}_o, \hat{\mathbf{y}}_o)$ is WR-efficient.

The set of all the WR-projections, denoted by $\Lambda_o$, is called the *projection set*.

## 2.2 The measurement of WR-RTS

From Banker et al. (2004), the WR-RTS generally has an unambiguous meaning for a WR-efficient DMU. In this case, the type and magnitude of the WR-RTS can be determined through the position(s) of the hyperplane(s) supporting the technology set at the DMU under evaluation. Under the uniqueness of the supporting hyperplane, obviously, the estimation of WR-RTS is straightforward. The difficulty arises, however, when the multiplicity occurs for the supporting hyperplanes. Naturally, the two potential questions arise here as to (1) what the cause is and (2) how to deal with it. In order to formally answer to these questions, we first introduce a new concept below.

**Definition 2.2.1** For a given WR-efficient activity $(\mathbf{x}, \mathbf{y})$, we define $(\mathbf{x}, \mathbf{y})$-*face* as the face of minimum dimension that contains $(\mathbf{x}, \mathbf{y})$ and denote it by $F_{(\mathbf{x}, \mathbf{y})}$.

It is worth noting that $(\mathbf{x}, \mathbf{y})$ is necessarily a relative interior point of $F_{(\mathbf{x},\mathbf{y})}$; since otherwise another face, namely $F$, will exist such that $(\mathbf{x}, \mathbf{y}) \in F \subseteq \partial\left(F_{(\mathbf{x},\mathbf{y})}\right) \subsetneqq F_{(\mathbf{x},\mathbf{y})}$, which is a contradiction by the definition. The following theorem demonstrates that the elements of $F_{(\mathbf{x},\mathbf{y})}$ are all WR-efficient.

**Theorem 2.2.1** $F_{(\mathbf{x},\mathbf{y})}$ is a strong face of $T_V^{WR}$.

*Proof.* Since $(\mathbf{x}, \mathbf{y})$ is WR-efficient, the *strictly complementary slackness condition* (Goldman & Tucker, 1956) in linear programming implies the existence of a strong supporting hyperplane of $T_V^{WR}$ that passes through $(\mathbf{x}, \mathbf{y})$. Let $H^S$ be such a hyperplane. We know that $(\mathbf{x}, \mathbf{y})$ is a relative interior point of $F_{(\mathbf{x},\mathbf{y})}$. Hence, Theorem 6.4 in Rockafellar (1970) implies that $H^S$ is binding over $F_{(\mathbf{x},\mathbf{y})}$ thereby the elements of $F_{(\mathbf{x},\mathbf{y})}$ are all WR-efficient. This completes the proof. ∎

In order to estimate the WR-RTS of a given WR-efficient DMU$_o$, we now consider the strong face $F_{(\mathbf{x}_o,\mathbf{y}_o)}$. It can be easily verified that $F_{(\mathbf{x}_o,\mathbf{y}_o)}$ is the intersection of $T_V^{WR}$ with all the supporting hyperplanes at $(\mathbf{x}_o, \mathbf{y}_o)$, i.e.,

$$F_{(\mathbf{x}_o,\mathbf{y}_o)} = T_V^{WR} \cap \bigcap_{H_o: \text{ supporting hyperplane at } (\mathbf{x}_o,\mathbf{y}_o)} H_o. \tag{5}$$

Hence, if the dimension of $F_{(\mathbf{x}_o,\mathbf{y}_o)}$ is $m+s-1$—i.e., it is a 'Full Dimensional Efficient Facet' (Olesen & Petersen, 1996, 2003), then a unique hyperplane is supporting at $(\mathbf{x}_o, \mathbf{y}_o)$. Otherwise, problem Type I arises, i.e., the multiplicity of supporting hyperplanes due to the non-full-dimensionality of $F_{(\mathbf{x}_o,\mathbf{y}_o)}$. To deal with problem Type I, the following theorem extends the method of Banker and Thrall (1992).

**Theorem 2.2.2** (Tone 2001; Korhonen et al., 2011) Let DMU$_o$ be WR-efficient and let $\bar{u}_0$ and $\underline{u}_0$ be the upper and lower bound of $u_0$, respectively, which are obtained by solving the following models:

$$\begin{aligned}
&\bar{u}_0(\underline{u}_0) = \max(\min) \ u_0 \\
&\text{subject to} \\
&\mathbf{u}^T \mathbf{y}_o - u_0 = 1, \\
&\mathbf{v}^T \mathbf{x}_o = 1, \\
&\mathbf{u}^T \mathbf{Y} - \mathbf{v}^T \mathbf{X} - u_0 \mathbf{1}_n \leq \mathbf{0}_n, \\
&\mathbf{v}^T \mathbf{P} \leq \mathbf{0}_k, \\
&\mathbf{u}^T \mathbf{Q} \leq \mathbf{0}_l, \\
&\mathbf{u} \geq \mathbf{0}_s, \ \mathbf{v} \geq \mathbf{0}_m, \ u_0 : \text{free in sign}.
\end{aligned} \tag{6}$$

Then,

- WR-IRS prevail at DMU$_o$ if and only if $\bar{u}_0 < 0$.

- WR-CRS prevail at DMU$_o$ if and only if $\underline{u}_0 \leq 0 \leq \bar{u}_0$.

- WR-DRS prevail at DMU$_o$ if and only if $\underline{u}_0 > 0$.

Having measured the WR-RTS for WR-efficient DMUs, the question arises now as to how to accomplish this task for the WR-inefficient ones. To determine the WR-RTS for a WR-inefficient DMU$_o$, it is first projected onto the efficient portion of the frontier as in (4). Then, the WR-RTS of its WR-projection is estimated as that of DMU$_o$.

While estimating the WR-RTS for a WR-inefficient DMU, one may encounter two types of problems. The first one is the occurrence of problem Type I for the given WR-projection $P_o$. This problem arises due to the non-full-dimensionality of $F_{P_o}$ and can be circumvented by Theorem 2.2.2. The second one—problem Type II—is the occurrence of multiple sets of supporting hyperplanes due to the occurrence of multiple WR-projections or, correspondingly, multiple $P_o$-faces. This problem is resolved if there exists a WR-projection $P_o^{\max}$ such that each supporting hyperplane at which passes through $\Lambda_o$ or, equivalently, that $F_{P_o} \subseteq F_{P_o^{\max}}$ for any $P_o \in \Lambda_o$; because under this condition a common set of supporting hyperplanes participate in the measurement of WR-RTS that are characterized by all the $P_o$-faces. In the immediately following section, we will prove that $P_o^{\max}$ always exists and will develop an LP model for its identification.

## 3 Our proposed approach

In this section, we develop an approach to cope effectively with the simultaneous occurrences of problems Types I and II. First, we introduce some concepts based on the optimal solution(s) of model (3).

**Definition 3.1** Let $\left(\theta^*, \boldsymbol{\lambda}^*, \boldsymbol{\pi}^*, \boldsymbol{\tau}^*, \mathbf{s}^{-*}, \mathbf{s}^{+*}\right)$ be an optimal solution to model (3) pertaining to a given WR-projection $P_o$. Then, the set of DMUs with positive $\lambda_j^*$ is defined as the *WR unary reference set (WR-URS)* for DMU$_o$ and is denoted by $R_{P_o}^U$ as

$$R_{P_o}^U = \left\{ \left(\mathbf{x}_j, \mathbf{y}_j\right) \middle| \lambda_j^* > 0 \right\}. \tag{7}$$

Each member of $R_{P_o}^U$ is referred to as a *reference DMU* of DMU$_o$. Tone (2001) proves that the reference units of DMU$_o$ are all WR-efficient in $T_V^{WR}$. Note that multiple WR-URSs occur if multiple optimal values take place for $\boldsymbol{\lambda}$. This requires us to define a reference set containing all the possible WR-URSs.

**Definition 3.2** The union of *all* the WR-URSs associated with a given projection $P_o$ is defined as the *WR maximal reference set (WR-MRS)* for DMU$_o$ and is denoted by $R_{P_o}^M$ as

$$R_{P_o}^M = \left\{ (\mathbf{x}_j, \mathbf{y}_j) \middle| \lambda_j^* > 0 \text{ in some optimal solution of (3) associated with } P_o \right\}. \tag{8}$$

Because of the occurrence of multiple WR-projections, we present the following definition to have a *unique* WR reference set containing all the possible reference DMUs.

**Definition 3.3** We define the union of *all* the WR-MRSs of DMU$_o$ as its *WR-global reference set (WR-GRS)* and denote it by $R_o^G$ as

$$R_o^G = \bigcup_{P_o \in \Lambda_o} R_{P_o}^M. \tag{9}$$

In the following theorem, we establish a relationship between Definitions 2.2.1 and 3.3 to demonstrate how problem Type II can be circumvented by identifying $R_o^G$.

**Theorem 3.1** There exists a WR-projection $P_o^{\max}$ such that $F_{P_o} \subseteq F_{P_o^{\max}}$ for any $P_o \in \Lambda_o$.

*Proof.* Let us define $\Omega_o = \left\{ \boldsymbol{\lambda} \middle| (\mathbf{X}\boldsymbol{\lambda}, \mathbf{Y}\boldsymbol{\lambda}) \in \Lambda_o \right\}$. Then, we have

$$R_o^G = \bigcup_{\boldsymbol{\lambda} \in \Omega_o} \left\{ \text{DMU}_j \middle| \lambda_j > 0 \right\}. \tag{10}$$

Since $\Omega_o$ is a non-empty set in $\mathbb{R}_{\geq 0}^n$, it has a *maximal element*, that is, an element with the maximum number of positive components. Let $\boldsymbol{\lambda}_o^{\max}$ be such an element. We now rewrite $\Omega_o$ equivalently as

$$\Omega_o = \left\{ \lambda \left| \begin{bmatrix} \mathbf{X} & -\mathbf{P} & \mathbf{0}_{m \times l} & \mathbf{I}_m & \mathbf{0}_{m \times s} \\ \mathbf{Y} & \mathbf{0}_{s \times k} & \mathbf{Q} & \mathbf{0}_{s \times m} & -\mathbf{I}_s \\ \mathbf{1}_n^T & \mathbf{0}_k^T & \mathbf{0}_l^T & \mathbf{0}_m^T & \mathbf{0}_s^T \\ \mathbf{0}_n^T & \mathbf{0}_k^T & \mathbf{0}_l^T & \mathbf{1}_m^T & \mathbf{1}_s^T \end{bmatrix} \begin{bmatrix} \lambda \\ \pi \\ \tau \\ \mathbf{s}^- \\ \mathbf{s}^+ \end{bmatrix} = \begin{bmatrix} \theta^* \mathbf{x}_o \\ \mathbf{y}_o \\ 1 \\ \mathbf{1}_m^T \mathbf{s}^{-*} + \mathbf{1}_s^T \mathbf{s}^{+*} \end{bmatrix}, \begin{pmatrix} \lambda \\ \pi \\ \tau \\ \mathbf{s}^- \\ \mathbf{s}^+ \end{pmatrix} \geq \begin{pmatrix} \mathbf{0}_n \\ \mathbf{0}_k \\ \mathbf{0}_l \\ \mathbf{0}_m \\ \mathbf{0}_s \end{pmatrix} \right\}. \quad (11)$$

From (11), $\Omega_o$ is a polyhedron and a convex set, accordingly. This implies that $\lambda_o^{\max}$ takes on positive values in any positive component of any $\lambda \in \Omega_o$. That is,

$$\left\{ j \mid \lambda_{oj}^{\max} > 0 \right\} = \bigcup_{\lambda \in \Omega_o} \left\{ j \mid \lambda_j > 0 \right\}. \quad (12)$$

Then, from (10) and (12), it holds that

$$R_o^G = \left\{ (\mathbf{x}_j, \mathbf{y}_j) \mid \lambda_{oj}^{\max} > 0 \right\}. \quad (13)$$

We now define $P_o^{\max}$ as

$$P_o^{\max} = \left( \mathbf{x}_o^{\max}, \mathbf{y}_o^{\max} \right) = \left( \mathbf{X} \lambda_o^{\max}, \mathbf{Y} \lambda_o^{\max} \right). \quad (14)$$

Then, $P_o^{\max} \in \Lambda_o$ by the definition of $\Omega_o$. From (13) and (14), $P_o^{\max}$ is a strict convex combination of the units in $R_o^G$. On the other hand, Theorem 2.2.1 implies that $F_{P_o^{\max}}$ is a face of $T_V^{WR}$ that contains $P_o^{\max}$. Therefore, by definition of face, we have that $conv\left( R_o^G \right) \subseteq F_{P_o^{\max}}$. Hence, the proof is complete by the fact that $\Lambda_o \subseteq conv\left( R_o^G \right)$. ∎

By the proof of Theorem 3.1, problem Type II can be resolved by identifying $\lambda_o^{\max}$. To make this identification, we propose the following LP problem:

$$\begin{aligned}&\max\ \mathbf{1}_n^T \boldsymbol{\mu}^1\\&\text{subject to}\end{aligned}$$

$$\begin{bmatrix} \mathbf{X} & -\mathbf{P} & \mathbf{0}_{m\times l} & \mathbf{I}_m & \mathbf{0}_{m\times s} \\ \mathbf{Y} & \mathbf{0}_{s\times k} & \mathbf{Q} & \mathbf{0}_{s\times m} & -\mathbf{I}_s \\ \mathbf{1}_n^T & \mathbf{0}_k^T & \mathbf{0}_l^T & \mathbf{0}_m^T & \mathbf{0}_s^T \\ \mathbf{0}_n^T & \mathbf{0}_k^T & \mathbf{0}_l^T & \mathbf{1}_m^T & \mathbf{1}_s^T \end{bmatrix} \begin{bmatrix} \boldsymbol{\mu}^1+\boldsymbol{\mu}^2 \\ \boldsymbol{\pi} \\ \boldsymbol{\tau} \\ \mathbf{s}^- \\ \mathbf{s}^+ \end{bmatrix} = \begin{bmatrix} \theta^*\mathbf{x}_o \\ \mathbf{y}_o \\ 1 \\ \mathbf{1}_m^T\mathbf{s}^{-*}+\mathbf{1}_s^T\mathbf{s}^{+*} \end{bmatrix}\eta, \quad (15)$$

$$\mathbf{0}_n \le \boldsymbol{\mu}^1 \le \mathbf{1}_n,\ \eta \ge 1,$$
$$\boldsymbol{\mu}^2 \ge \mathbf{0}_n,\ \boldsymbol{\pi} \ge \mathbf{0}_k,\ \boldsymbol{\tau} \ge \mathbf{0}_l,\ \mathbf{s}^- \ge \mathbf{0}_m,\ \mathbf{s}^+ \ge \mathbf{0}_s.$$

**Theorem 3.2** Let $\left(\boldsymbol{\mu}^{1*},\boldsymbol{\mu}^{2*},\boldsymbol{\pi}^*,\boldsymbol{\tau}^*,\mathbf{s}^{-*},\mathbf{s}^{+*},\eta^*\right)$ be an optimal solution to model (15). Then, $\lambda_o^{\max} = \dfrac{1}{\eta^*}\left(\boldsymbol{\mu}^{1*}+\boldsymbol{\mu}^{2*}\right).$

*Proof.* It is obvious that $\dfrac{1}{\eta^*}\left(\boldsymbol{\mu}^{1*}+\boldsymbol{\mu}^{2*}\right)\in \Omega_o$. Without loss of generality, let only the first $p$ components of this vector are positive. We then claim that $\lambda_o^{\max} = \dfrac{1}{\eta^*}\left(\boldsymbol{\mu}^{1*}+\boldsymbol{\mu}^{2*}\right)$. By the convexity of $\Omega_o$, this is equivalent to say that $\dfrac{1}{\eta^*}\left(\boldsymbol{\mu}^{1*}+\boldsymbol{\mu}^{2*}\right)$ takes positive value in any positive component of any $\lambda \in \Omega_o$. Assume by contradiction that this is not true, that is, there exists a $\bar{\lambda}\in \Omega_o$ for which $\bar{\lambda}_j > 0$ for some $j \in \{p+1,...,n\}$. Then, we have

$$\begin{bmatrix} \mathbf{X} & -\mathbf{P} & \mathbf{0}_{m\times l} & \mathbf{I}_m & \mathbf{0}_{m\times s} \\ \mathbf{Y} & \mathbf{0}_{s\times k} & \mathbf{Q} & \mathbf{0}_{s\times m} & -\mathbf{I}_s \\ \mathbf{1}_n^T & \mathbf{0}_k^T & \mathbf{0}_l^T & \mathbf{0}_m^T & \mathbf{0}_s^T \\ \mathbf{0}_n^T & \mathbf{0}_k^T & \mathbf{0}_l^T & \mathbf{1}_m^T & \mathbf{1}_s^T \end{bmatrix} \begin{bmatrix} \boldsymbol{\mu}^{1*}+\boldsymbol{\mu}^{2*}+\bar{\boldsymbol{\lambda}} \\ \boldsymbol{\pi}^*+\bar{\boldsymbol{\pi}} \\ \boldsymbol{\tau}^*+\bar{\boldsymbol{\tau}} \\ \mathbf{s}^{-*}+\bar{\mathbf{s}}^- \\ \mathbf{s}^{+*}+\bar{\mathbf{s}}^+ \end{bmatrix} = \begin{bmatrix} \theta^*\mathbf{x}_o \\ \mathbf{y}_o \\ 1 \\ \mathbf{1}_m^T\mathbf{s}^{-*}+\mathbf{1}_s^T\mathbf{s}^{+*} \end{bmatrix}(\eta^*+1). \quad (16)$$

We define $\boldsymbol{\pi}' := \boldsymbol{\pi}^* + \overline{\boldsymbol{\pi}}$, $\boldsymbol{\tau}' := \boldsymbol{\tau}^* + \overline{\boldsymbol{\tau}}$, $\mathbf{s}^{-\prime} := \mathbf{s}^{-*} + \overline{\mathbf{s}}^-$, $\mathbf{s}^{+\prime} := \mathbf{s}^{+*} + \overline{\mathbf{s}}^+$, $\eta' := \eta^* + 1$ and

$$\mu_j^{1\prime} := \begin{cases} \mu_j^{1*}, & j = 1,...,p, \\ \min\{1, \overline{\lambda}_j\}, & j = p+1,...,n, \end{cases}, \quad \mu_j^{2\prime} := \begin{cases} \mu_j^{2*}, & j = 1,...,p, \\ \overline{\lambda}_j - \mu_j^{1\prime}, & j = p+1,...,n. \end{cases} \tag{17}$$

Then, $\left(\boldsymbol{\mu}^{1\prime}, \boldsymbol{\mu}^{2\prime}, \boldsymbol{\pi}', \boldsymbol{\tau}', \mathbf{s}^{-\prime}, \mathbf{s}^{+\prime}, \eta'\right)$ is a feasible solution to (15) whose objective function value is strictly greater than $\mathbf{1}_n^T \boldsymbol{\mu}^{1*}$. This contradicts the optimality of $\left(\boldsymbol{\mu}^{1*}, \boldsymbol{\mu}^{2*}, \boldsymbol{\pi}^*, \boldsymbol{\tau}^*, \mathbf{s}^{-*}, \mathbf{s}^{+*}, \eta^*\right)$ and proves our claim. ∎

Having identified $\boldsymbol{\lambda}_o^{\max}$, one can also identify $R_o^G$ by (13).

We now propose the following two-stage procedure for measuring the WR-RTS of a WR-inefficient DMU so as to deal with the simultaneous occurrence of problems Types I and II:

**Stage 1** Solve model (15) and yield $P_o^{\max}$ by (14).

**Stage 2** Measure the WR-RTS of $P_o^{\max}$ by using Theorem 2.2.2.

In summary, the WR-RTSs of all DMUs can be estimated by the execution of the following three main steps:

**Step 1** We first evaluate each DMU via model (3) to obtain the efficiency score and the sum of the optimal input and output slacks. Based on the obtained results for all the observed DMUs, we then partition them into two disjoint groups. The first group includes WR-efficient and WR-inefficient DMUs with zero sums of slacks. Each unit in this group has a unique WR-projection and, hence, it may face only with problem Type I. The second group also contains the remaining WR-inefficient

DMUs, i.e., those with non-zero sums of slacks. For each unit in this group, multiple WR-projections may occur and, consequently, it may suffer from either of problems Types I or II or both.

**Step 2** We measure the WR-RTS of DMUs in the first group by applying Theorem 2.2.2 to, the DMU under evaluation if it is WR-efficient or, to its unique WR-projection if it is WR-inefficient.

**Step 3** We use our proposed two-stage procedure to estimate the WR-RTS of the units in the second group.

## 4 An empirical application

To demonstrate the ready applicability of our proposed method, we conduct an illustrative empirical analysis based on a data set of 80 Iranian secondary schools, denoted by S1–S80, participating in the Third International Mathematics and Science Study (TIMSS) supported by UNESCO. To carry out all the computations, we have developed a computer program using the GAMS optimization software that can be found in Appendix A.

The data set, as used and analyzed by Korhonen et al. (2011), consists of four inputs: the income level of the parents ($x_1$), educational facilities of the students at home ($x_2$), contribution of the parents towards the school programs ($x_3$) and the education level of the teachers ($x_4$); and one output, i.e., the school's GPA based on individual students' grades of each school in the TIMSS study ($y$). The input–output data can be found in Table 2 in Korhonen et al. (2011). As in Korhonen et al. (2011), we assume that $k = 3$ and the matrix $P$ is given as

$$\left.\begin{array}{r} v_1 - \frac{1}{3}v_2 \leq 0 \\ v_2 - \frac{1}{2}v_3 \leq 0 \\ v_3 - v_4 \leq 0 \end{array}\right\} \quad \Rightarrow \quad P = \begin{pmatrix} 1 & 0 & 0 \\ -\frac{1}{3} & 1 & 0 \\ 0 & -\frac{1}{2} & 1 \\ 0 & 0 & -1 \end{pmatrix}_{4\times 3}.$$

**Step 1:** *Evaluation of schools via model (3)*

We have evaluated all the schools via the WR model (3) and have presented the results in Table 1. As can be seen, schools S5, S9, S26, S56, S57, S77 and S80 are all rated as WR-efficient and the remaining schools are WR-inefficient.

Out of 75 WR-inefficient schools, 30 schools have zero sums of slacks. Therefore, the set of all schools is partitioned into two groups. The first group consists of 37 schools (7 WR-efficient and 30 WR-inefficient), where each faces only with problem Type I. The second group contains 43 WR-inefficient schools that may suffer from either of problems Types I or II or both. In all the tables, we have colored the labels of schools in the first and second groups in green and red, respectively.

**Step 2:** *Measuring the WR-RTSs of schools in the first group*

As previously mentioned, schools in the first group are WR-efficient or WR-inefficient and their WR-RTS can be measured in accordance with Theorem 2.2.2. For WR-efficient schools, however, we have solved both minimization and maximization forms of model (6) to obtain the lower and upper bounds of $u_0$, respectively. The obtained bounds together with the WR-RTSs of schools are all provided in Table 2. According to Theorem 2.2.2, schools S5, S26, S56 and S80 exhibit constant WR-RTS, schools S9 and S57 increasing WR-RTS and school S77 decreasing WR-RTS.

Table 1 The results obtained from model (3) for 80 schools

| DMU | $\theta^*$ | $\mathbf{1}_m^T \mathbf{s}^{-*} + \mathbf{1}_s^T \mathbf{s}^{+*}$ | DMU | $\theta^*$ | $\mathbf{1}_m^T \mathbf{s}^{-*} + \mathbf{1}_s^T \mathbf{s}^{+*}$ |
|---|---|---|---|---|---|
| S1 | 0.6667 | 82.5622 | S41 | 0.8988 | 0 |
| S2 | 0.6505 | 77.1072 | S42 | 0.7253 | 9.6181 |
| S3 | 0.5973 | 5.7104 | S43 | 0.642 | 0 |
| S4 | 0.5785 | 66.2883 | S44 | 0.6886 | 0.0504 |
| S5 | 1 | 0 | S45 | 0.791 | 0 |
| S6 | 0.6635 | 13.001 | S46 | 0.6454 | 0 |
| S7 | 0.8061 | 98.439 | S47 | 0.7445 | 1.0731 |
| S8 | 0.6444 | 87.2777 | S48 | 0.7 | 0 |
| S9 | 1 | 0 | S49 | 0.6649 | 0 |
| S10 | 0.759 | 54.9407 | S50 | 0.7188 | 8.0192 |
| S11 | 0.7484 | 0 | S51 | 0.6883 | 0 |
| S12 | 0.9327 | 0 | S52 | 0.6676 | 58.5 |
| S13 | 0.684 | 0 | S53 | 0.6037 | 18.5478 |
| S14 | 0.631 | 0 | S54 | 0.771 | 13.7189 |
| S15 | 0.6314 | 9.4 | S55 | 0.6376 | 57.8916 |
| S16 | 0.7012 | 65.7235 | S56 | 1 | 0 |
| S17 | 0.7748 | 0 | S57 | 1 | 0 |
| S18 | 0.6473 | 46.2338 | S58 | 0.8077 | 0 |
| S19 | 0.6916 | 30 | S59 | 0.6624 | 0 |
| S20 | 0.6611 | 35.7 | S60 | 0.8331 | 0 |
| S21 | 0.6705 | 24.547 | S61 | 0.796 | 0 |
| S22 | 0.723 | 12.2596 | S62 | 0.8821 | 63.4544 |
| S23 | 0.8543 | 53.8568 | S63 | 0.7768 | 0 |
| S24 | 0.7539 | 0 | S64 | 0.7057 | 0 |
| S25 | 0.8021 | 6.5305 | S65 | 0.6048 | 64.5015 |
| S26 | 1 | 0 | S66 | 0.6418 | 61.5521 |
| S27 | 0.6602 | 12.9995 | S67 | 0.6739 | 0 |
| S28 | 0.6586 | 42.2138 | S68 | 0.6513 | 2.6223 |
| S29 | 0.9414 | 0 | S69 | 0.8 | 7.4087 |
| S30 | 0.7674 | 33.7742 | S70 | 0.6936 | 0 |
| S31 | 0.7184 | 53.3647 | S71 | 0.7718 | 10.9831 |
| S32 | 0.7685 | 29.5648 | S72 | 0.7996 | 0 |
| S33 | 0.6981 | 49.0425 | S73 | 0.7421 | 39.1132 |
| S34 | 0.8074 | 3.7996 | S74 | 0.6598 | 31.068 |
| S35 | 0.8037 | 0 | S75 | 0.7824 | 0 |
| S36 | 0.8143 | 0 | S76 | 0.9992 | 4.0789 |
| S37 | 0.6428 | 0 | S77 | 1 | 0 |
| S38 | 0.7342 | 7.8291 | S78 | 0.9985 | 0 |
| S39 | 0.6952 | 34.2956 | S79 | 0.8777 | 0 |

| | | | | | | |
|---|---|---|---|---|---|---|
| S40 | 0.817 | 0 | | S80 | 1 | 0 |

Table 2 The WR-RTSs for 7 WR-efficient schools in the first group

| DMU | $\underline{u}_0$ | $\overline{u}_0$ | WR-RTS |
|---|---|---|---|
| S5 | -1 | 0.1699 | C |
| S9 | -1 | -0.4588 | I |
| S26 | -1 | 0.2623 | C |
| S56 | -1 | 2.4818 | C |
| S57 | -1 | -0.0388 | I |
| S77 | 0.5735 | 1E+10 | D |
| S80 | -0.1169 | 4.9070 | C |

Note: C: Constant WR-RTS; I: Increasing WR-RTS; and D: Decreasing WR-RTS

To estimate the WR-RTS of each WR-inefficient school, we have first determined its unique WR-projection from the optimal solution of model (3). Then, we have determined its WR-RTS by computing the lower and upper bounds of $u_0$. The results obtained for 30 WR-inefficient schools are all presented in Table 3. The results reveal that, out of these 30 schools, schools S75, S78 and S79 operate under constant WR-RTS, and the remaining ones under increasing WR-RTS.

Table 3 The WR-RTSs for 30 WR-inefficient schools in the first group

| DMU | $\hat{x}_1$ | $\hat{x}_2$ | $\hat{x}_3$ | $\hat{x}_4$ | $\hat{y}$ | $\underline{u}_0$ | $\overline{u}_0$ | WR-RTS |
|---|---|---|---|---|---|---|---|---|
| S11 | 5.4867 | 13.2084 | 3.8324 | 4.4905 | 409.3 | -1 | -0.4484 | I |
| S12 | 7.4814 | 13.9999 | 4.5153 | 3.7307 | 419.1 | -1 | -0.4425 | I |
| S13 | 8.0103 | 12.6862 | 4.4627 | 4.1039 | 419.3 | -1 | -0.4424 | I |
| S14 | 5.4921 | 12.9233 | 3.5693 | 4.8614 | 415.7 | -1 | -0.4445 | I |
| S17 | 6.8498 | 12.6459 | 3.9118 | 4.6487 | 421 | -1 | -0.4414 | I |
| S24 | 9.6787 | 11.8251 | 4.3992 | 4.5236 | 437.5 | -1 | -0.4319 | I |
| S29 | 10.8733 | 12.7631 | 4.9861 | 3.7657 | 440.8 | -1 | -0.4301 | I |
| S35 | 9.508 | 11.6522 | 4.1627 | 4.8223 | 441.5 | -1 | -0.4297 | I |
| S36 | 11.358 | 14.1296 | 3.5003 | 4.8855 | 466.6 | -0.7705 | -0.0744 | I |
| S37 | 9.4173 | 12.2036 | 3.7067 | 5.1426 | 449.3 | -0.7771 | -0.077 | I |
| S40 | 9.8242 | 12.1109 | 3.9761 | 4.9017 | 447.3 | -0.7779 | -0.0774 | I |
| S41 | 10.166 | 11.6901 | 4.4891 | 4.4938 | 440.2 | -1 | -0.4304 | I |
| S43 | 5.6917 | 12.7569 | 3.5028 | 4.9945 | 419.3 | -1 | -0.4424 | I |
| S45 | 7.9557 | 12.4824 | 3.955 | 4.7459 | 433 | -0.5328 | -0.4345 | I |
| S46 | 6.981 | 12.2101 | 3.5641 | 5.1635 | 430.7 | -1 | -0.4358 | I |
| S48 | 11.2669 | 12.1873 | 4.676 | 4.2002 | 444.8 | -0.7788 | -0.0778 | I |
| S49 | 5.7323 | 13.5455 | 4.2179 | 3.9894 | 402.7 | -1 | -0.4524 | I |
| S51 | 8.4896 | 12.5114 | 4.5119 | 4.1297 | 422.9 | -1 | -0.4403 | I |
| S58 | 10.0008 | 12.2002 | 4.0114 | 4.8459 | 447.9 | -0.7776 | -0.0773 | I |
| S59 | 7.1441 | 12.0538 | 3.4908 | 5.299 | 434.1 | -1 | -0.4339 | I |
| S60 | 5.31 | 13.1977 | 3.6547 | 4.6761 | 411.6 | -0.8169 | -0.4657 | I |
| S61 | 10.4519 | 12.7273 | 4.7758 | 3.9798 | 441.6 | -0.5279 | -0.4297 | I |
| S63 | 11.479 | 13.6957 | 3.8363 | 4.6608 | 461.4 | -0.7725 | -0.0752 | I |
| S64 | 7.5954 | 12.7212 | 4.3116 | 4.2342 | 418.8 | -1 | -0.4427 | I |
| S67 | 11.111 | 11.728 | 4.9425 | 4.0434 | 438.7 | -1 | -0.4313 | I |
| S70 | 9.0588 | 12.3026 | 4.5694 | 4.1616 | 427.2 | -1 | -0.4378 | I |
| S72 | 10.3384 | 12.5255 | 3.9816 | 4.7976 | 450.8 | -0.7765 | -0.0768 | I |
| S75 | 10.9182 | 13.5616 | 3.3681 | 6 | 494.8 | -0.1338 | 0.1508 | C |
| S78 | 14.5748 | 15.589 | 3.7302 | 6 | 542.3 | -0.1235 | 0.1358 | C |
| S79 | 13.0312 | 15.4629 | 3.3088 | 5.4682 | 506 | -0.1312 | 0.0319 | C |

Note: C: Constant WR-RTS; I: Increasing WR-RTS; and D: Decreasing WR-RTS

**Step 3:** *Measuring the WR-RTSs of schools in the second group*

For each WR-inefficient school in the second group, we have applied our proposed tow-stage procedure discussed in Section 3. We have first obtained the maximal intensity vector $\lambda_o^{\max}$ by solving model (15). Table 4 exhibits the WR-GRSs together with their associated intensity vectors. Except school S77, all the remaining WR-efficient schools appear in the WR-GRSs of WR-inefficient schools. For example, schools S5, S9 and S77 with the respective weights of 0.3334, 0.1126 and 0.554 constitute the WR-GRS of school S1. This means that target inputs and outputs of school S1 must be adjusted as a convex combination of inputs and outputs of S5, S9 and S77, in order it becomes WR-efficient.

Table 4 The WR-GRSs for 43 schools in the second group

| DMU | WR-global reference set | | | | | |
|---|---|---|---|---|---|---|
| | S5 | S9 | S26 | S56 | S57 | S80 |
| S1 | 0.3334 | 0.1126 | 0 | 0 | 0.554 | 0 |
| S2 | 0.6021 | 0.2818 | 0 | 0 | 0.1161 | 0 |
| S3 | 0.0335 | 0.9665 | 0 | 0 | 0 | 0 |
| S4 | 0.3141 | 0.139 | 0 | 0 | 0.5469 | 0 |
| S6 | 0.281 | 0.719 | 0 | 0 | 0 | 0 |
| S7 | 0 | 0.2245 | 0 | 0.7755 | 0 | 0 |
| S8 | 0 | 0.4489 | 0 | 0.1339 | 0.4173 | 0 |
| S10 | 0.2769 | 0.6086 | 0 | 0 | 0.1146 | 0 |
| S15 | 0 | 1 | 0 | 0 | 0 | 0 |
| S16 | 0.3997 | 0.6003 | 0 | 0 | 0 | 0 |
| S18 | 0.9129 | 0 | 0 | 0 | 0.0871 | 0 |
| S19 | 0 | 1 | 0 | 0 | 0 | 0 |
| S20 | 0 | 1 | 0 | 0 | 0 | 0 |
| S21 | 0.9916 | 0.0084 | 0 | 0 | 0 | 0 |
| S22 | 1 | 0 | 0 | 0 | 0 | 0 |
| S23 | 0.5427 | 0.4573 | 0 | 0 | 0 | 0 |
| S25 | 0.4063 | 0.4644 | 0 | 0 | 0.1293 | 0 |
| S27 | 0 | 0.961 | 0 | 0.039 | 0 | 0 |
| S28 | 0.6346 | 0 | 0 | 0 | 0.3654 | 0 |
| S30 | 0 | 0.0696 | 0 | 0.9304 | 0 | 0 |
| S31 | 0.8735 | 0 | 0 | 0 | 0.1265 | 0 |
| S32 | 0.2138 | 0.7862 | 0 | 0 | 0 | 0 |
| S33 | 0.0944 | 0.1997 | 0 | 0 | 0.7058 | 0 |
| S34 | 0.0725 | 0 | 0.6993 | 0 | 0.2282 | 0 |
| S38 | 0.2026 | 0.7404 | 0 | 0 | 0.057 | 0 |
| S39 | 0.7806 | 0 | 0 | 0 | 0.2194 | 0 |
| S42 | 0.9012 | 0 | 0 | 0 | 0.0988 | 0 |
| S44 | 0.8316 | 0 | 0 | 0 | 0 | 0.1684 |
| S47 | 0.7947 | 0 | 0.0443 | 0 | 0 | 0.161 |
| S50 | 0.6258 | 0.3742 | 0 | 0 | 0 | 0 |
| S52 | 0 | 1 | 0 | 0 | 0 | 0 |
| S53 | 0.2301 | 0.7699 | 0 | 0 | 0 | 0 |
| S54 | 0.3129 | 0.6616 | 0 | 0 | 0.0255 | 0 |
| S55 | 0.449 | 0.551 | 0 | 0 | 0 | 0 |
| S62 | 0 | 0 | 0 | 0.8239 | 0 | 0.1761 |
| S65 | 0.4192 | 0.0632 | 0 | 0 | 0.5176 | 0 |
| S66 | 0.3115 | 0.6885 | 0 | 0 | 0 | 0 |
| S68 | 0.5975 | 0.4025 | 0 | 0 | 0 | 0 |

| DMU | | | | | | |
|---|---|---|---|---|---|---|
| S69 | 0.4001 | 0.4703 | 0 | 0 | 0.1296 | 0 |
| S71 | 0.3154 | 0 | 0 | 0 | 0.6846 | 0 |
| S73 | 0.1832 | 0.8168 | 0 | 0 | 0 | 0 |
| S74 | 0 | 0.9515 | 0 | 0.0415 | 0.0071 | 0 |
| S76 | 0 | 0 | 0.5037 | 0 | 0 | 0.4963 |

After identifying the intensity vector $\lambda_o^{max}$ and the WR-GRS for each school, we have then determined the corresponding WR-projection $P_o^{max}$ as in (14). The results are given in the second to sixth columns of Table 5. Finally, we have determined the lower and upper bounds of $u_0$ by using model (6) and have estimated the WR-RTS based on these bounds. The obtained bounds together with the WR-RTS of hospitals are provided in the last three columns of Table 5. As can be observed from Table 5, out of 43 WR-inefficient hospitals, five hospitals exhibit constant WR-RTS and the remaining hospitals increasing WR-RTS.

**Table 5** The WR-RTSs for 43 schools in the second group

| DMU | $x_1^{max}$ | $x_2^{max}$ | $x_3^{max}$ | $x_4^{max}$ | $y_1^{max}$ | $\underline{u}_0$ | $\bar{u}_0$ | WR-RTS |
|---|---|---|---|---|---|---|---|---|
| S1 | 9.1335 | 11.8924 | 4.2206 | 4.6667 | 436.3622 | -1 | -0.4326 | I |
| S2 | 6.8306 | 12.2277 | 3.514 | 5.2042 | 430.4072 | -1 | -0.436 | I |
| S3 | 4.3004 | 13.9163 | 3.9665 | 4.0669 | 394.2104 | -1 | -0.4577 | I |
| S4 | 9.0252 | 11.9569 | 4.2327 | 4.6283 | 434.9883 | -1 | -0.4334 | I |
| S6 | 5.0429 | 13.2976 | 3.719 | 4.562 | 407.601 | -1 | -0.4494 | I |
| S7 | 8.6204 | 15.6286 | 4.7755 | 3.2245 | 436.139 | -1 | -0.4327 | I |
| S8 | 7.9256 | 13.3214 | 4.5511 | 3.8661 | 419.4777 | -1 | -0.4423 | I |
| S10 | 5.844 | 13.0443 | 3.8377 | 4.5538 | 412.7407 | -1 | -0.4463 | I |
| S15 | 4.2 | 14 | 4 | 4 | 392.4 | -1 | -0.4588 | I |
| S16 | 5.3991 | 13.0008 | 3.6003 | 4.7994 | 414.0235 | -1 | -0.4455 | I |
| S18 | 7.5571 | 11.5174 | 3.1742 | 5.8258 | 445.8642 | -1 | -0.0776 | I |
| S19 | 4.2 | 14 | 4 | 4 | 392.4 | -1 | -0.4588 | I |
| S20 | 4.2 | 14 | 4 | 4 | 392.4 | -1 | -0.4588 | I |
| S21 | 7.1749 | 11.5209 | 3.0084 | 5.9833 | 446.047 | -1 | -0.4272 | I |
| S22 | 7.2 | 11.5 | 3 | 6 | 446.5 | -1 | 0.1699 | C |
| S23 | 5.8281 | 12.6432 | 3.4573 | 5.0854 | 421.7608 | -1 | -0.441 | I |
| **S25** | 6.3366 | 12.6869 | 3.7229 | 4.8126 | 420.4305 | -1 | -0.4417 | I |
| **S27** | 4.4223 | 14.0819 | 4.039 | 3.961 | 394.5995 | -1 | -0.4574 | I |
| **S28** | 8.6982 | 11.5731 | 3.7309 | 5.2691 | 443.8324 | -1 | -0.0779 | I |
| **S30** | 9.5032 | 15.9538 | 4.9304 | 3.0696 | 444.8742 | -1 | -0.4278 | I |
| **S31** | 7.7186 | 11.5253 | 3.253 | 5.747 | 445.5766 | -1 | -0.0776 | I |
| **S32** | 4.8413 | 13.4656 | 3.7862 | 4.4275 | 403.9648 | -1 | -0.4516 | I |
| **S33** | 9.4948 | 12.1405 | 4.6114 | 4.1889 | 430.5425 | -1 | -0.4359 | I |
| **S34** | 11.7022 | 14.4828 | 3.4564 | 4.8443 | 469.8 | -0.7693 | -0.0739 | I |
| **S38** | 5.2128 | 13.3624 | 3.8545 | 4.4051 | 406.0291 | -1 | -0.4503 | I |
| **S39** | 8.0995 | 11.5439 | 3.4388 | 5.5612 | 444.8984 | -1 | -0.0777 | I |
| **S42** | 7.6051 | 11.5198 | 3.1976 | 5.8024 | 445.7788 | -1 | -0.0776 | I |
| **S44** | 8.9013 | 12.4433 | 3.1684 | 6 | 468.6 | -0.1402 | 0.1606 | C |
| **S47** | 9.0516 | 12.5874 | 3.161 | 5.9557 | 469.2 | -0.1401 | 0.0344 | C |
| **S50** | 6.0774 | 12.4355 | 3.3742 | 5.2516 | 426.2554 | -1 | -0.4383 | I |
| **S52** | 4.2 | 14 | 4 | 4 | 392.4 | -1 | -0.4588 | I |
| **S53** | 4.8903 | 13.4248 | 3.7699 | 4.4602 | 404.8478 | -1 | -0.4511 | I |
| **S54** | 5.3196 | 13.1592 | 3.7126 | 4.6257 | 410.5189 | -1 | -0.4476 | I |
| **S55** | 5.547 | 12.8775 | 3.551 | 4.898 | 416.6916 | -1 | -0.4439 | I |
| **S62** | 11.2032 | 16.2761 | 4.8239 | 3.5283 | 471.5 | -0.0991 | 2.1102 | C |
| **S65** | 9.1325 | 11.7616 | 4.0984 | 4.8384 | 439.3015 | -1 | -0.4309 | I |
| **S66** | 5.1345 | 13.2213 | 3.6885 | 4.623 | 409.2521 | -1 | -0.4484 | I |
| **S68** | 5.9924 | 12.5064 | 3.4025 | 5.1949 | 424.7223 | -1 | -0.4392 | I |
| **S69** | 6.3202 | 12.7018 | 3.7295 | 4.8002 | 420.1087 | -1 | -0.4419 | I |
| **S71** | 10.0069 | 11.6369 | 4.3692 | 4.6308 | 441.5023 | -1 | -0.0783 | I |
| **S73** | 4.7497 | 13.5419 | 3.8168 | 4.3665 | 402.3132 | -1 | -0.4526 | I |
| **S74** | 4.4864 | 14.0708 | 4.0485 | 3.9585 | 395.068 | -1 | -0.4571 | I |
| **S76** | 14.7817 | 16.3949 | 3.4963 | 5.4963 | 529.6 | -0.1261 | 0.2333 | C |

Note: C: Constant WR-RTS; I: and Increasing WR-RTS

# 5 Concluding remarks

Since weight restrictions are necessary for real-life applications of DEA, the problem of measuring returns to scale in this framework (WR-RTS) is certainly of great importance. As is well known, the major difficulty underlying the problem of WR-RTS measurement comes from the occurrence of multiple supporting hyperplanes. This current study contributes to the literature by resolving this difficulty. For a WR-efficient point *P*, the concept of *P*-face is first introduced as a face of minimum dimension that contains *P*. This concept is then employed to exactly determine two potential sources of the origin of the multiplicity issue—(i) the non-full-dimensionality of the *P*-face (problem Type I) and (ii) the occurrence of multiple *P*-faces (problem Type II). While problem Type I occurs for both WR-efficient and WR-inefficient DMUs, problem Type II arises only for WR-inefficient DMUs.

To cope with the simultaneous occurrences of problems Types I and II, three types of WR-reference set—WR-unary reference set (WR-URS), WR-maximal reference set (WR-MRS) and WR-global reference set (WR-GRS), are introduced. The WR-GRS is a unique reference set that is defined as the union of WR-MRSs associated with all the WR-projections. For a given WR-projection *P*, the WR-MRS is the union of the WR-URSs associated with *P* and each WR-URS contains WR-efficient DMUs that are active in a specific convex combination producing *P*. The reason for introducing three types of WR-reference set is to have a mathematically well-defined definition of WR-reference set. Indeed, by this introduction, a clear distinction is made between the uniquely-found WR-reference set—the WR-GRS—and the two types of WR-reference set for which multiplicity may occur—the WR-URS and the WR-MRS.

Finally, a linkage between the *P*-face and WR-GRS notions is established wherein it is proved that if $P^{max}$ is a strict convex combination of the units in the WR-GRS, then $P^{max}$-face contains all the possible *P*-faces. This implies that problem Type II is circumvented if $P^{max}$ is used for the measurement of the WR-RTS. Then, the problem under consideration is reduced to identifying $P^{max}$ or, equivalently, to identifying a maximal element of the set of intensity vectors associated with the optimal solutions of the WR-model. An LP problem is developed to find such an element to overcome problem Type II. By applying Tone's (2001) method to $P^{max}$, problem Type I is then resolved.

To test the full efficiency of a DMU and identify its reference units with production trade-offs, Podinovski (2007) developed a three-stage procedure. The third stage of his procedure involves the identification of reference set by means of an LP model. Since his proposed LP model does not identify all the possible reference units, we suggest the use of our proposed method in his procedure.